\documentclass{amsart}
\usepackage{amssymb,latexsym}
\numberwithin{equation}{section}
\newtheorem{theorem}{Theorem}[section]
\newtheorem{proposition}[theorem]{Proposition}

\newtheorem{remark}[theorem]{Remark}

\begin{document}

\pagenumbering{arabic}
\pagestyle{headings}
\def\sof{\hfill\rule{2mm}{2mm}}
\def\ls{\leq}
\def\gs{\geq}
\def\SS{\mathcal S}
\def\qq{{\bold q}}
\def\qe{{\bold e}}
\def\qf{{\bold f}}
\def\txx{{\frac1{2\sqrt{x}}}}
\def\mn{\mbox{-}}

\title{ {\sc continued fractions, statistics, and generalized patterns}}
   
\author{Toufik Mansour}
\maketitle
\begin{center}{LABRI, Universit\'e Bordeaux I,\\ 
              351 cours de la Lib\'eration, 33405 Talence Cedex\\
		{\tt toufik@labri.fr} }
\end{center}
%
\section*{Abstract}
Recently, Babson and Steingrimsson (see \cite{BS}) introduced generalized permutations 
patterns that allow the requirement that two adjacent letters in a 
pattern must be adjacent in the permutation. 

Following \cite{BCS}, let $e_k\pi$ (respectively; $f_k\pi$) be the number of 
the occurrences of the generalized pattern $12\mn3\mn\dots\mn k$ 
(respectively; $21\mn3\mn\dots\mn k$) in $\pi$. In the present note, 
we study the distribution of the statistics $e_k\pi$ and $f_k\pi$ 
in a permutation avoiding the classical pattern $1\mn3\mn2$.

Also we present an applications, which relates the Narayana numbers, Catalan numbers, and 
increasing subsequences, to permutations avoiding the classical pattern $1\mn3\mn2$ 
according to a given statistics on $e_k\pi$, or on $f_k\pi$. 
\section{Introduction}
{\bf Permutation patterns:} let $\alpha\in S_n$ and $\tau\in S_k$ be two permutations. We say that $\alpha$ 
{\em contains} $\tau$ if there exists a subsequence $1\leq i_1<i_2<\cdots<i_k\leq n$ 
such that $(\alpha_{i_1},\dots,\alpha_{i_k})$ is order-isomorphic to $\tau$; 
in such a context $\tau$ is usually called a {\em pattern} or ({\em classical pattern}). 
We say $\alpha$ {\em avoids}
$\tau$, or is $\tau$-{\em avoiding}, if such a subsequence does not exist. 
The set of all $\tau$-avoiding permutations in $S_n$ is denoted $S_n(\tau)$. 
For an arbitrary finite collection of patterns $T$, we say that $\alpha$
avoids $T$ if $\alpha$ avoids any $\tau\in T$; the corresponding subset of 
$S_n$ is denoted $S_n(T)$.

While the case of permutations avoiding a single pattern has attracted
much attention, the case of multiple pattern avoidance remains less
investigated. In particular, it is natural, as the next step, to consider
permutations avoiding pairs of patterns $\tau_1$, $\tau_2$. This problem
was solved completely for $\tau_1,\tau_2\in S_3$ (see \cite{SS}), for 
$\tau_1\in S_3$ and $\tau_2\in S_4$ (see \cite{W}), and for 
$\tau_1,\tau_2\in S_4$ (see \cite{B1,Km} and references therein).   
Several recent papers \cite{CW,MV,Kr,MV3,MV2} deal with the case 
$\tau_1\in S_3$, $\tau_2\in S_k$ for various pairs $\tau_1,\tau_2$. Another
natural question is to study permutations avoiding $\tau_1$ and containing
$\tau_2$ exactly $t$ times. Such a problem for certain $\tau_1,\tau_2\in S_3$ 
and $t=1$ was investigated in \cite{R}, and for certain $\tau_1\in S_3$, 
$\tau_2\in S_k$ in \cite{RWZ,MV,Kr,MV3}. The tools involved in these papers 
include continued fractions, Chebyshev polynomials, and Dyck paths.\\

{\bf Generalized permutation patterns:}
in \cite{BS} Babson and Steingrimsson introduced generalized permutation 
patterns that allow the requirement that two adjacent letters in a pattern
must be adjacent in the permutation. This idea in 
introducing these patterns was study of Mahonian statistics. 

We write a classical pattern with dashes between any two adjacent letters of the pattern, 
say $1324$, as $1\mn3\mn2\mn4$, and if we write, say $24\mn1\mn3$, then we mean that 
if this pattern occurs in permutation $\pi$, then the letters in the 
permutation $\pi$ that correspond to $2$ and $4$ are adjacent (see \cite{C}).
For example, the permutation $\pi=35421$ has only two occurrences of the pattern $23\mn1$, 
namely the subsequnces $352$ and $351$, whereas $\pi$ has four occurrences of the 
pattern $2\mn3\mn1$, namely the subsequnces $352$, $351$, $342$ and $341$. 

Claesson \cite{C} gave a complete answer for the number of permutations avoiding any 
single $3$-leters generalized pattern with exactly one adjacent pair of letters. 
Later, Claesson and Mansour \cite{CM}
presented a complete solution for the number of permutations avoiding any double 
$3$-letters generalized patterns with exactly one adjacent pair of letters. 
Besides, Kitaev \cite{Ki} investigate simultaneous avoidance of two or more $3$-letters 
generalized patterns without internal dashes.\\

On the other hand, Robertson, Wilf and Zeilberger \cite{RWZ} showed a simple continued 
fraction that records the joint distribution of the patterns $1\mn2$ and $1\mn2\mn3$ 
on $1\mn3\mn2$-avoiding permutations. Recently, generalization of this theorem given, 
by Mansour and Vainshtein \cite{MV}, by Krattenthaler \cite{Kr}, 
by Jani and Rieper \cite{JR}, and by Br\"and\"en, Claesson and Steingrimsson \cite{BCS}.

Mansour \cite[Th. 2.1 and Th. 2.9]{M} 
presented an analog of this theorem (\cite{MV}) by replace 
generalized patterns with classical patterns. \\

In the present note, 
we generalize \cite[Th. 2.1]{M} and \cite[Th. 2.9]{M} and 
give an analog for \cite{BCS} by replace generalized patterns with classical patterns. 
This by use the same argument proof of the main results in \cite{BCS} with a 
simple changes.  In the last section, we present an applications for 
our results.

\section{Main results}
For all $k\geq 1$, we denote by $a_k(\pi)$ the number of the occurrences 
of the pattern $1\mn2\mn3\mn\dots\mn k$ in $\pi$. In \cite[Th. 1]{BCS} proved 
the following.

\begin{theorem} {\rm (P. Br\"and\'en, A. Claesson, and E. Steingrimsson, \cite[Th. 1]{BCS})}
\label{bcs1}
The $\sum\limits_{\pi\in\mathcal{S}(1\mn3\mn2)} \prod\limits_{k\geq 1} x_k^{a_k(\pi)}$ 
is given by the following continued fraction:
$$
\dfrac{1}
{1-\dfrac{x_1^{{0\choose 0}}}
{1-\dfrac{x_1^{{1\choose 0}}x_2^{{1\choose 1}}}
{1-\dfrac{x_1^{{2\choose 0}}x_2^{{2\choose 1}}x_3^{{2\choose 2}}}
{\ddots}}}}
$$
in which the $(n+1)$st numerator is $\prod\limits_{k=0}^n x_{k+1}^{{n\choose k}}$.
\end{theorem}

For all $k\geq 3$, we denote by $e_k(\pi)$ (respectively; $f_k(\pi)$) 
the number of the occurrences of the generalized pattern $12\mn3\mn\dots\mn k$ 
(respectively; $21\mn3\mn\dots\mn k$) in $\pi$. Besides, $e_2(\pi)$, $f_2(\pi)$, and 
$e_1(\pi)=f_1(\pi)$ are denote the number of occurrences of the pattern 
$12$, $21$ and $1$, respectively, and extended that by $e_0(\pi)=f_0(\pi)=1$.

Now let us give an analog for Theorem \ref{bcs1} by used generalized patterns with 
simply changes in the proof of Theorem \ref{bcs1}.

\begin{theorem}
\label{cc}
The $\sum\limits_{\pi\in\mathcal{S}(1\mn3\mn2)} \prod\limits_{k\geq 1} x_k^{e_k(\pi)}$ is given by the 
following continued fraction:
$$
\dfrac{1}
{1-x_1+x_1x_2^{{0\choose 0}}-\dfrac{x_1x_2^{{0\choose 0}}}
{1-x_1+x_1x_2^{{1\choose 0}}x_3^{{1\choose 1}}-\dfrac{x_1x_2^{{1\choose 0}}x_3^{{1\choose 1}}}
{1-x_1+x_1x_2^{{2\choose 0}}x_3^{{1\choose 1}}x_4^{{2\choose 2}}-\dfrac{x_1x_2^{{2\choose 0}}x_3^{{2\choose 1}}x_4^{{2\choose 2}}}
{\ddots}}}}
$$
in which the $(n+1)$st numerator is $x_1\cdot\prod\limits_{k=0}^n x_{k+2}^{{n\choose k}}$.
\end{theorem}
\begin{proof}
Let $\pi$ any nonempty permutation avoiding $1\mn2\mn3$ such that $\pi=(\pi',n,\pi'')$ and 
$j=\pi^{-1}(n)$; so every element in $\pi'$ is greater than every element in $\pi''$. 
Thus $\pi',\pi''\in\mathcal{S}(1\mn3\mn2)$, so  
$$e_k(\pi)=\left\{ \begin{array}{ll} 
			e_k(\pi')+e_{k-1}(\pi')+e_k(\pi''),\quad & \mbox{for all}\ k\geq 3;\\
			e_2(\pi')+e_2(\pi'')+\delta_{\pi',\emptyset}, & k=2;\\
			e_1(\pi')+e_1(\pi'')+1,			 & k=1
	           \end{array}\right.,$$
where $\delta_{\pi',\emptyset}$ is the Kronecker delta.

It follows that the generating function 
	$$C(x_1,x_2,\dots)=\sum_{\pi\in\mathcal{S}(1\mn3\mn2)} \prod_{k\geq 1} x_k^{e_k(\pi)}$$
satisfies (see \cite[Th. 2.1]{M})
	$$C(x_1,x_2,\dots)=1+x_1C(x_1,x_2,\dots)+x_1x_2C(x_1,x_2,\dots) (C(x_1,x_2x_3,x_3x_4,\dots)-1),$$
which means that ,
	$$C(x_1,x_2,x_3,x_4,\dots)=\frac{1}{1-x_1+x_1x_2-x_1x_2C(x_1,x_2x_3,x_3x_4,\dots)}$$
and the theorem follows by induction.
\end{proof}

Similarly, it is easy to see 
$$f_k(\pi)=\left\{ \begin{array}{ll} 
			f_k(\pi')+f_{k-1}(\pi')+f_k(\pi''),\quad 	& \mbox{for all}\ k\geq 3;\\
			f_2(\pi')+f_2(\pi'')+\delta_{\pi'',\emptyset}, 	& k=2;\\
			f_1(\pi')+f_1(\pi'')+1,			 	& k=1
	           \end{array}\right.,$$
where $\delta_{\pi'',\emptyset}$ is the Kronecker delta. So by using 
\cite[Th. 2.9]{M} and the argument proof of Theorem \ref{bcs1} 
we get the following.

\begin{theorem}
\label{dd}
The $\sum\limits_{\pi\in\mathcal{S}(1\mn3\mn2)} \prod\limits_{k\geq 1} x_k^{f_k(\pi)}$ is given by the 
following continued fraction:
$$
1-\dfrac{x_1}
{x_1x_2^{{0\choose 0}}-\dfrac{1}
{1-\dfrac{x_1}
{x_1x_2^{{1\choose 0}}x_3^{{1\choose 1}}-\dfrac{1}
{1-\dfrac{x_1}
{x_1x_2^{{2\choose 0}}x_3^{{2\choose 1}}x_4^{{2\choose 2}}-\dfrac{1}{\ddots}}}}}}
.$$
\end{theorem}

Following \cite{BCS}, we define $\mathcal{A}$ be 
the ring of all infinite matrices with finite number of non zero 
entries in each row, that is, 
$$\mathcal{A}=\{A:\mathcal{N}\times\mathcal{N}\rightarrow\mathcal{Z}|\forall n
	(A(n,k)=0\ \mbox{for almost every\ } k)\},$$
with multiplication defined by $(AB)(n,k)=\sum_{i\geq 1} A(n,i)B(i,k)$.
With each $A\in\mathcal{A}$ we now associate a family of statistics $\{<\qq,A_k>\}_{k\geq 1}$ 
defined on $\mathcal{S}(1\mn3\mn2)$, where $\qq=(q_1,q_2,q_3,\dots)$ and 
	$$<\qq,A_k>=\sum_{i\geq 1} A(i,k)q_i.$$

Following \cite{BCS}, let us define mathematical objects with respect $A$ as follows. 
Let $\qq=(q_1,q_2,\dots)$ where the $q_i$ are indeterminates; for each 
$A\in\mathcal{A}$ and $\pi\in\mathcal{S}(1\mn3\mn2)$ we define three objects as follows:

The {\em weight} $\eta(\pi,A;\qq)$, the {\em weight} $\mu(\pi,A;\qq)$,  
and the {\em weight} $\nu(\pi,A;\qq)$ of $\pi$ 
	with respect $A$, by 
	$$\eta(\pi,A;\qq)=\prod_{k\geq 1} q_{k}^{<{\bf a},A_k>\pi},\quad
	  \mu(\pi,A;\qq)=\prod_{k\geq 1} q_{k}^{<\qe,A_k>\pi},\quad 
	  \nu(\pi,A;\qq)=\prod_{k\geq 1} q_{k}^{<\qf,A_k>\pi},$$
	respectively, where ${\bf a}=(a_1,a_2,\dots)$, $\qe=(e_1,e_2,\dots)$ and $\qf=(f_1,f_2,\dots)$. 

The {\em generating function with respect $A$} of the, statistics $\{<{\bf a},A_k>\}_{k\geq 1}$, 
statistics $\{<\qe,A_k>\}_{k\geq 1}$, and statistics $\{<\qf,A_k>\}_{k\geq 1}$ by 
	$$F_A(\qq)=\sum_{\pi\in\mathcal{S}(1\mn3\mn2)} \eta(\pi,A;\qq),\quad
	  G_A(\qq)=\sum_{\pi\in\mathcal{S}(1\mn3\mn2)} \mu(\pi,A;\qq),\quad
	  H_A(\qq)=\sum_{\pi\in\mathcal{S}(1\mn3\mn2)} \nu(\pi,A;\qq),$$
respectively. 

The {\em continued fractions with respect $A$}, by
	$$C_A(\qq)=\dfrac{1}
	{1-\dfrac{\prod\limits_{k\geq 0} q_{k+1}^{A(1,k)}}
	{1-\dfrac{\prod\limits_{k\geq 0} q_{k+1}^{A(2,k)}}
	{\ddots}}},$$

	$$D_A(\qq)=\dfrac{1}
	{1-q_1+q_1\prod\limits_{k\geq 1} q_{k+1}^{A(1,k)}-\dfrac{q_1\prod\limits_{k\geq 1} q_{k+1}^{A(1,k)}}
	{1-q_1+q_1\prod\limits_{k\geq 1} q_{k+1}^{A(2,k)}-\dfrac{q_1\prod\limits_{k\geq 1} q_{k+1}^{A(2,k)}}
	{\ddots}}},$$
	and by
	$$E_A(\qq)=1-\dfrac{q_1}
	{q_1\prod_{k\geq 1} q_{k+1}^{A(1,k)}-\dfrac{1}
	{1-\dfrac{q_1}
	{q_1\prod_{k\geq 1} q_{k+1}^{A(2,k)}-\dfrac{1}
	{\ddots}}}}.$$

As a remark (see \cite[p. 3]{BCS}): definition of the ring $\mathcal{A}$,  
and the fact that $a_i(\pi)=e_i(\pi)=f_i(\pi)=0$ for all $i=|\pi|+1,|\pi|+2,\dots$ 
yields the product in the weight is finite. 

The second step in \cite{BCS}, proved the following. 

\begin{theorem} {\rm (P. Br\"and\'en, A. Claesson, and E. Steingrimsson, \cite[Th. 2]{BCS})}
\label{bcs2}
For $A\in\mathcal{A}$, $F_A(\qq)=C_{BA}(\qq)$ where $B=\left[\binom{i}{j} \right]$, and 
conversely $C_A(\qq)=F_{B^{-1}A}(\qq)$. 
\end{theorem}

By above definitions with using the proof of Theorem \ref{bcs2} 
we get an analog of Theorem \ref{bcs2} for the statistics $e_k$ and $f_k$ 
as follows. 

\begin{theorem}
\label{mm}
Let $A\in\mathcal{A}$; then  
	$$G_A(\qq)=D_{BA}(\qq),\quad\quad H_A(\qq)=E_{BA}(\qq),$$
	$$D_A(\qq)=G_{B^{-1}A}(\qq),\quad\quad E_A(\qq)=H_{B^{-1}A}(\qq),$$
where $B=\left[ {i\choose j}\right]$, $\binom{n}{k}=0$ for all $k>0$ or $n>k$. 
\end{theorem}
\begin{proof}
By using definitions we get that (here, we using the same argument proof 
in \cite[Th. 2]{BCS})
$$\begin{array}{ll}
\mu(\pi,A;\qq)	&=\prod_{k\geq 1} q_k^{<\qe,A_k>\pi}\\
		&=\prod_{k\geq 1}\prod_{j\geq 1} q_k^{A(j,k)e_j(\pi)}\\
		&=\prod_{j\geq 1}\left( \prod_{k\geq 1} q_k^{(A(j,k)} \right)^{e_j(\pi)}.
\end{array}$$
Let $x_{j+1}=\prod_{k\geq 1} q_k^{A(j,k)}$; Theorem \ref{cc} yields 
	$$\prod_{j\geq 1} x_{j+1}^{{n-1\choose {j-1}}}=
	  \prod_{j\geq 1} \left( \prod_{k\geq 1} q_k^{A(j,k)} \right)^{{{n-1}\choose{j-1}}}=
	  \prod_{k\geq 1} q_k^{\left( {{n-1}\choose 0}, {{n-1}\choose 1}, {{n-1}\choose 2},\dots\right),A_k>}.$$
so, again, by definitions $G_A(\qq)=D_{BA}(\qq)$. Observing that 
$B^{-1}=\left[ (-1)^{i-j}{i\choose j} \right]\in\mathcal{A}$ we also 
obtain $D_A(\qq)=G_{B^{-1}A}(\qq)$. 

Similarly, we have $E_{BA}(\qq)=H_A(\qq)$ and $E_A(\qq)=H_{B^{-1}A}(\qq)$. 
\end{proof}

\begin{remark}
The general approach which described in \cite{BCS} on the statistics $a_k(\pi)$ 
its work with others statistics. Its work on the statistics 
$e_k(\pi)$ and on the statistics $f_k(\pi)$. So the natural question to ask is 
the following: If there any descriptions for all the statistics such 
this approach will work? 
Here we failed to give an answer for this.
\end{remark}
\section{Application}
In the current section,  we present an examples for application of Theorem \ref{mm}. 
Some of these examples are related known 
continued fraction to the statistics $e_k$ or $f_k$, but others relate these statistics 
to others combinatorial objects.

\subsection{Narayana numbers} 
Let $N(n,k)=\frac{1}{n}{n\choose k}{n\choose {k+1}}$ be the Narayana numbers, 
and the corresponding generating function for the Narayana numbers we denote by 
$N(x,t)$. Then 
	$$N(x,t):=\sum_{n,k\geq 0} N(n,k)x^kt^n=1+xtN^2(x,t)-xtN(x,t)+tN(x,t).$$
This allow us to express $N(x,t)$ as a continued fraction:
	$$N(x,t)=\dfrac{1}{1-\dfrac{t}{1-\dfrac{tx}{1-\dfrac{t}{1-\dfrac{tx}{\ddots}}}}}.$$

\begin{proposition}
The number $N(n,k)$ equals the number of permutations $\pi\in S_n(132)$ with $e_2(\pi)=k$.
\end{proposition}
\begin{proof}
Let $A(n,k)=\delta_{(n,k),(1,1)}+\delta_{(n,k),(2,2)}$ where $\delta$ is Kronecker delta, so
by applying Theorem \ref{mm} we get that  
$$N'(x,t):=\sum_{\pi\in\mathcal{S}(1\mn3\mn2)} x^{e_2(\pi)}t^{|\pi|}=\frac{1}{1-t-xt+xtN'(x,t)}=
\dfrac{1}{1-t-xt+\dfrac{xt}{1-t-xt+\dfrac{xt}{\ddots}}},$$
so $N'(x,t)$ satisfies the same functional equation as $N(x,t)$, hence $N'(x,t)=N(x,t)$.
\end{proof}

Again, the same argument work also to statistics on $f_k$ as follows.

\begin{proposition}
The number $N(n,k)$ equals the number of permutations $\pi\in S_n(132)$ with $f_2(\pi)=k$.
\end{proposition}
\begin{proof}
Let $A(n,k)=\delta_{(n,k),(1,1)}+\delta_{(n,k),(2,2)}$ where $\delta$ is Kronecker delta, so
by applying Theorem \ref{mm} we get that 
$$N''(x,t):=\sum_{\pi\in\mathcal{S}(1\mn3\mn2)} x^{f_2(\pi)}t^{|\pi|}=1-\dfrac{t}{xt-\dfrac{1}{1-tN''(x,t)}}=
1-\dfrac{t}{xt-\dfrac{1}{1-\dfrac{t}{xt-\dfrac{1}{1-\dfrac{t}{\ddots}}}}},$$
so $N''(x,t)$ satisfies the same functional equation as $N(x,t)$, hence $N''(x,t)=N(x,t)$.
\end{proof}

\subsection{Increasing subsequences}. 
As consequence of \cite{BCS}, we define as follows. 
The subsequence  $\pi_{i_1}\pi_{i_1+1}\pi_{i_2}\dots \pi_{i_k}$ ($k\geq 2$)
of $\pi$ is called {\em $2$-increasing subsequence} if $\pi_{i_j}<\pi_{i_{j+1}}$, $i_j<i_{j+1}$ 
and $i_1+1<i_2$. Hence, the total number of $2$-increasing subsequences in a permutation 
is counting by $e_2+e_3+\dots$. An application of Theorem \ref{mm}  
gives the following continued fraction 
for the distribution of $e_2+e_3+\dots$:
$$\sum_{\pi\in\mathcal{S}(1\mn3\mn2)} x^{e_2(\pi)+e_3(\pi)+\dots}t^{|\pi|}=
\dfrac{1}{1-t(1-x)-\dfrac{xt}{1-t(1-x^2)-\dfrac{x^2t}{1-t(1-x^4)-\dfrac{x^4t}{\ddots}}}}.$$

The subsequence  $\pi_{i_1}\pi_{i_1+1}\pi_{i_2}\dots \pi_{i_k}$ ($k\geq 2$)
of $\pi$ is called {\em almost $2$-increasing subsequence} if $\pi_{i_j}<\pi_{i_{j+1}}$, $i_j<i_{j+1}$ for $j=2,3,\dots,k-1$, 
$i_1+1<i_2$ and $\pi_{i_1+1}<p_{i_1}$. Hence, the total number of almost $2$-increasing subsequences in a permutation 
is counting by $f_2+f_3+\dots$. An application of Theorem \ref{mm} gives the following continued fraction 
for the distribution of $f_2+f_3+\dots$:
$$\sum_{\pi\in\mathcal{S}(1\mn3\mn2)} x^{f_2(\pi)+f_3(\pi)+\dots}t^{|\pi|}=
1-\dfrac{t}{xt-\dfrac{1}{1-\dfrac{t}{x^2t-\dfrac{1}{1-\dfrac{t}{x^4t-\dfrac{1}{1-\dfrac{t}{x^8t-\ddots}}}}}}}$$

\subsection{Catalan numbers}  The $n$th Catalan number is given by $C_n=\frac{1}{n+1}{{2n}\choose n}$, 
and the corresponding generating function is given by $C(x)=\frac{1-\sqrt{1-4x}}{2x}$. Theorem 
\ref{mm} yields for the statistic $s=0$ ($s(\pi)=0$ for all $\pi$) the following.
The generating function $C(x)$ for the number of permutations avoiding $1\mn3\mn2$ can be 
expressed, again, in terms of continued fractions:
$$C'(x):=\dfrac{1}{1-\dfrac{x}{1-\dfrac{x}{1-\dfrac{x}{\ddots}}}},\quad\quad 
C''(x):=1-\dfrac{x}{x-\dfrac{1}{1-\dfrac{x}{x-\dfrac{1}{\ddots}}}}$$
In the above two cases $C'(x)=C''(x)=C(x)$.

\end{document}